\documentclass{article}[12pt]
\usepackage[hmargin=2.5 cm]{geometry}

\usepackage{amsmath}
\usepackage{amssymb}
\usepackage{amsthm}
\usepackage{enumerate}
\usepackage{color}
\usepackage{clrscode}
\usepackage{multirow}
\usepackage{url}

\usepackage{tikz}
\usepackage{wrapfig}
\usetikzlibrary{arrows}
\tikzstyle{block}=[draw opacity=0.7,line width=1.4cm]
\usepackage{xspace}

\newcommand{\ran}{Randi\'{c} index\xspace}
\newcommand{\rann}{Randi\'{c} index\xspace}
\newcommand{\er}{Erd\H{o}s-R\'{e}nyi }
\newcommand{\ern}{Erd\H{o}s-R\'{e}nyi}
\newcommand{\p}{+}
\newcommand{\m}{-}

\newcommand{\Z}{\mathbb{Z}}

\newcommand{\R}{\mathbb{R}}

\newcommand{\cP}{\mathcal{P}}

\newcommand{\pf}{{\it Proof. }}

\newtheorem{definition}{Definition}[section]
\newtheorem{theorem}[definition]{Theorem}
\newtheorem{corollary}[definition]{Corollary}

\begin{document}
%%%%%%%%%%%%%%%%

\setcounter{page}{1} %

\title{Algorithm and Complexity for a Network Assortativity Measure}
\author{Sarah J. Kunkler\thanks{Department of Mathematics, College of William and Mary, Williamsburg, VA 23187-8795, USA, \texttt{sjkunkler@email.wm.edu}} \and M. Drew LaMar\thanks{Department of Biology, College of William and Mary, Williamsburg, VA 23187-8795, USA, \texttt{mdlama@wm.edu}} \and Rex K. Kincaid\thanks{Department of Mathematics, College of William and Mary, Williamsburg, VA 23187-8795, USA, \texttt{rrkinc@wm.edu}} \and David Phillips\thanks{Department of Mathematics, United States Naval Academy, Annapolis, MD 21401, USA, \texttt{dphillip@usna.edu}}}
\date{\today}

\maketitle

\abstract{ We show that finding a graph realization with the minimum
  \ran for a given degree sequence is solvable in polynomial time by
  formulating the problem
  as a minimum weight perfect $b$-matching problem.  However, the
  realization found via this reduction is not guaranteed to be
  connected.  
  Approximating the minimum weight $b$-matching problem
  subject to a connectivity constraint is shown to be NP-Hard. 
  For instances in which the optimal solution to the minimum \ran
  problem is not connected, we 
  describe a heuristic to connect the graph using pairwise edge
  exchanges that preserves the degree sequence.  In our computational
  experiments, the heuristic performs well and the \ran of the
  realization after our heuristic is within $3$\% of the unconstrained
  optimal value on average. Although we focus on minimizing the \ran, our 	  
  results extend to maximizing the \ran as well. Applications of the \ran to synchronization of neuronal networks controlling respiration in mammals and to
  normalizing cortical thickness networks in diagnosing individuals with
  dementia are provided.}

%%%%%%%%%%%%%%%%%%%%%%%%%%%%%%%%%%%%%%%%%%%%%%%%%%%%%%%%%
\section{Introduction}%1
%%%%%%%%%%%%%%%%%%%%%%%%%%%%%%%%%%%%%%%%%%%%%%%%%%%%%%%%%
Networks are pervasive in the sciences.  For example, they are used
in ecology to represent food webs and in engineering and computer
science to design high quality internet router connections.  Depending
on the application, one particular graph property may be more
important than another. Oftentimes, a desired property is to have a
connected graph or to optimize a particular metric while constrained
to connected graphs \cite{Kincaid11}.

One of these measures, the \ran of a graph, developed by Milan
Randi\'{c}, was originally used in chemistry \cite{randic75}. The \ran
of a graph can be thought of as an assortativity measure. A network is
described as disassortative if high-degree nodes are predominantly
attached to low-degree nodes \cite{li2005towards}. Minimizing the
\rann, in many instances, will produce a graph with disassortativity
\cite{li2005towards}.  Why is this optimization problem of interest?
Li et al. \cite{li2005towards}, in the design of an internet router
network, found that networks that maximized throughput also had small
values for the \ran.  In addition, the \ran has been shown
\cite{Kincaid09} to correlate with synchronization, an important
property in many network applications.  We explore this correlation
for directed graphs in section \ref{sec:netsync}. Optimizing the \ran is also useful in analyzing imaging data of cortical thicknesses which we discuss in section \ref{sec:ct}. Our focus is to investigate
algorithms that minimize the \ran of a graph over all connected
realizations while keeping the degrees of the nodes fixed. 
However, our complexity results and algorithms also extend 
to the case of maximizing the \ran as well.

\subsection{Notation and Definitions}
\label{sec:notation}

We assume the reader to have a knowledge of graph theory (see, e.g.,
\cite{West}). We let $\R$ and $\Z$ denote the set of real numbers and integers, respectively. We consider an undirected graph, $G =
(N,E)$, which consists of nodes, $N$, and edges, $E$. We assume that our graph is simple, i.e., there are no self-loops and no multi-edges.
The degree of a node is defined as $d_i(G) := |\{j:(i,j) \in E\}|$.
We denote the node-node adjacency matrix by $A(G)$. The \textit{degree sequence} is the list of the degrees of all the nodes in a graph, which 
we represent as  $d(G)=(d_1(G),d_2(G),...,d_n(G))$.
A sequence of non-negative integers is considered 
\textit{graphic} if it is the degree sequence of a graph.  
Degree sequences can correspond to more than one adjacency matrix or graph. We call these graphs different 
realizations of the degree sequence. When the particular graph is clear
from context, we omit $G$ in the previous definitions.

Let nodes $u,v \in N$. We say that $u$ and $v$ are \textit{connected} if there exists a path from $u$ to $v$. 
A graph is connected if for all $u\in N$ there exists a path to every other node.

\vskip .1 true in
\begin{definition}
The \ran of a graph $G=(N,E)$ is defined as $$R_\alpha(G)=\sum_{(v_i,v_j)\in E} (d_i \cdot d_j)^{\alpha},$$ where $\alpha \in \R - \{ 0 \}$.
\end{definition}
\vskip .1 true in
{\noindent}
\!\!We consider the case when $\alpha = 1$ which has received considerable attention in a variety of contexts~\cite{gutman1973graph,li2005towards,beichl2008measuring}. For convenience, we define $R(G) = R_1(G)$.
A natural optimization problem is: \\

{\noindent}\textbf{Minimum Randi\'{c} Index Problem.} Given a graphic degree sequence  
what is a graph realization with the minimum \rann with $\alpha=1$? \\

We define the {\it connected minimum \ran problem} as the minimum \ran
problem with the additional constraint of minimizing over all connected realizations. 

\vskip .1 true in
\begin{definition}
\label{def:pbmatch}
For a graph $G = (N,E)$ and a positive integer vector $b=(b_1,\cdots,b_n) \in
\Z^n$, a {\it perfect $b$-matching} is a subset of edges $M
\subseteq E$ such that for node $v_i \in N$, the degree of $v_i$ in the
graph $(N,M)$ is $b_i$.
\end{definition}
\vskip .1 true in

\noindent
An associated optimization problem is: \\

{\noindent}\textbf{Minimum Weight Perfect b-Matching Problem.}  Given a positive integer
vector $b$, a graph $G=(N,E)$ and a set of edge weights $w:E \rightarrow \R$, find a perfect
$b$-matching with minimum weight. 
\vskip .1 true in
In section 3, we will see that the minimum \ran problem is equivalent to the
minimum weight perfect $b$-matching problem on a complete graph $G$ with an appropriate choice of weights.  We will also show that, using an arbitrary graph $G$, by constraining the minimum weight perfect $b$-matching problem to connected perfect $b$-matchings, the problem becomes NP-Hard.

\subsection{Network Measures of Assortativity}
\label{sec:randic}

The \ran of a graph was originally defined in chemistry. 
In 1975, the chemist Milan
Randi\'{c} \cite{randic75} proposed the index $R_{\alpha}(N)$
for the cases $\alpha = -1 \hbox{ and } \alpha = -1/2$ under the name {\it branching index}. He explained the utility of $R$ in measuring the extent of branching of
the carbon-atom skeleton of saturated hydrocarbons.  
His approach ``reveals some inherent relationships between [structures] which can be traced to connectivity" \cite{randic75}. It is sometimes 
referred to as the connectivity index by scholars in chemistry~\cite{hansen2003variable}.
Bollob\'{a}s and Erd\"{o}s~\cite{Erdos} generalized this index by allowing $\alpha$ to take
on any non-zero real number.  A survey of results for the Randi\'{c} index can be found in \cite{LSW08}. 

The \ran with $\alpha=1$ is important enough to have been ``discovered'' several times in the literature.
For example, the second Zagreb index, defined by Gutman~\cite{gutman1973graph} and also used in chemistry, is actually $R(G)$.  In 2005, Li et al. \cite{li2005towards} investigated what they called the s-metric of a graph which, seemingly unknown to them, is actually $R(G)$. They used $R(G)$ to differentiate between graph realizations of a given 
degree sequence following a power law distribution \cite{li2005towards} in the design of internet router networks. 
For a fixed degree distribution they plotted $R(G)$ versus throughput for hundreds of graph realizations.  They noted that $R(G)$ measures the ``hub-like core" of a graph and is maximized when high-degree nodes are connected to other high degree nodes (assortative).  Conversely, the minimum values of $R(G)$ were predominantly associated with networks that maximized throughput and were dissassortative. In 2008, Beichl and Cloteaux investigated how well random networks generated with a chosen $R(G)$ can model the structure of real networks such as the 
Internet. The graphs produced by optimizing $R(G)$ resulted in better models than the ones that used simple uniform sampling \cite{beichl2008measuring}. 

\subsection{Random Graph Classes}
\label{sec:graphs}

Our computational experiments require random graphs. 
We make use of three types of graphs: \ern, geometric and scale-free. The structure 
of these graphs depends on the parameters chosen. \\

{\noindent}\textbf{\er Graphs.}
A number of nodes $n$ and a probability of connection $p$ are chosen. A uniform random number on the interval $[0,1]$ is generated for each possible edge. If the 
number generated for an edge is less than $p$ then the edge is added. \\

{\noindent}\textbf{Geometric Graphs.}
A number of nodes $n$ and a radius $r$ is chosen.  Each node $v_i$ is placed uniformly at random in the unit square, giving coordinates $(x_{i},y_{i})$.  We connect nodes $v_i$ and $v_j$ if $(x_i - x_j)^2 + (y_i -y_j)^2 \leq r^2$ \cite{taylor2009contest}.  \\

{\noindent}\textbf{Scale-Free Graphs.}
A preferential attachment algorithm is used to create graphs whose degree sequences
follow a power-law distribution. Following the convention in the literature we will refer
to these graphs as ``scale-free''.  A number of nodes $n$ is chosen. New nodes are added and 
connected to existing nodes, based on a probability proportional to the current degree of the nodes, 
until you reach $n$ nodes, making it more likely that a new node will be connected to a higher degree node \cite{taylor2009contest}.
The algorithm allows a minimum node degree to be specified. \\

\section {Formulation and Complexity}

In this section, we formulate %reduce
the minimum \ran problem as a \textit{minimum weight perfect
b-matching problem}, which is solvable in polynomial time
\cite{schrijver2003combinatorial}. Note that this problem does not
enforce connectivity. We then show that approximating the minimum \ran
problem with connectivity is NP-Hard.

\subsection{The $b$-Matching Problem }

Consider a graph $G = (N,E)$, a positive integer vector $b=(b_1,\cdots,b_n) \in
\Z^n$ and $M \subseteq E$, a perfect b-matching.  For a given $b$-matching, 
$M$, the graph induced by $M$ is $(N,M)$. We denote the set of perfect $b$-matchings of
a graph $G$ by $\cP_b(G)$.  For edge weights $w:E \rightarrow \R$, the
{\it minimum weight perfect $b$-matching problem} is finding the perfect
$b$-matching with minimum weight, i.e., to calculate
\begin{equation}
  \label{eq:perfectb}
  M^\ast(G) := \arg\min\left\{\sum_{e \in M} w(e) : M \in \cP_b(G)\right\}.
\end{equation}
For example, let G be the undirected, weighted graph below
\begin{center}
\begin{tikzpicture}[node distance=20mm]
        \tikzstyle{every node}=%
        [%
          fill=blue!50!black!20,%
          draw=blue!50!black,%
          minimum size=8mm,%
          circle,%
          thick,%
          inner sep=0pt%
        ]

        \node (1)              {$v_1$};
        \node (2) [right of=1] {$v_2$};
        \node (3) [below of=1] {$v_3$};
        \node (4) [below of=2] {$v_4$};

        \path [thick] [every node/.style={font=\sffamily\small}]    
        		      (1) edge node[above] {3} (2) 
                         (1) edge node [xshift= -0.7 cm, yshift = 0.4 cm] {2} (4) 
                         (3) edge node [xshift= 0.7 cm, yshift = 0.4 cm] {7} (2) 
                         (2) edge node[right] {4} (4) 
                         (3) edge node [below]{1}(4);
      \end{tikzpicture} \\
\end{center}
and let $b = (2,1,1,2)$ for nodes $v_1$, $v_2$, $v_3$ and $v_4$ respectively. 
We select $b_i$ edges that will connect to the $i$th node and that
will produce the minimum weight. Therefore, the matching induces the graph $G'$ below.
\begin{center}
\begin{tikzpicture}[node distance=20mm]
        \tikzstyle{every node}=%
        [%
          fill=blue!50!black!20,%
          draw=blue!50!black,%
          minimum size=8mm,%
          circle,%
          thick,%
          inner sep=0pt%
        ]

        \node (1)              {$v_1$};
        \node (2) [right of=1] {$v_2$};
        \node (3) [below of=1] {$v_3$};
        \node (4) [below of=2] {$v_4$};

        \path [thick][every node/.style={font=\sffamily\small}]    
        		      (1) edge node[above] {3} (2) 
                         (1) edge node [xshift= -0.35 cm] {2} (4) 
                         (3) edge node [below]{1}(4);

      \end{tikzpicture} \\
\end{center}
Note that for this example the solution $G'$ is the only perfect $b$-matching for $G$.  Let an instance of the minimum \ran problem be given with a positive integer vector $b
\in \Z^n$ and graph $G=(N,E)$.  Then the graphs induced
by the matchings are feasible subgraph realizations for the minimum
\ran problem. Also, the edge sets
of the feasible subgraph realizations of the minimum \ran problem are
perfect $b$-matchings. Thus, the set of feasible $b$-matchings on $G$ is identical to
the set of feasible subgraph realizations on $G$ to the minimum \ran
problem on $G$. 
Thus, to formulate an instance of minimum \ran problem as a minimum weight perfect
$b$-matching problem, set
\begin{equation}
  \label{eq:w}
  w_{ij} = b_i \cdot b_j.
\end{equation}
Therefore, we can create an instance of a minimum weight perfect
$b$-matching to solve the minimum \ran problem. Since the $b$-matching
problem can be solved in polynomial time, finding the minimum \ran of
a graph can also be done in polynomial time.  Note that this method
does not enforce connectivity. 

We show that even approximating the connected minimum \ran problem is
NP-Hard. We first define approximation algorithms (see \cite{Vazirani01} for
further details about approximation algorithms). Let $S \subset \R^n$
and $f:S \rightarrow \R$ be a given feasibility set and objective
function, respectively. Define an $\alpha$-approximation algorithm for
the minimization problem $v^* = \min_{x \in S} f(x)$ as a polynomial
time algorithm that finds a solution $y \in S$ with $f(y) \leq \alpha
v^*$. We say that we can {\it approximate} a minimization problem if
there exists an $\alpha$ such that an $\alpha$-approximation algorithm
exists. Note that $\alpha \geq 1$ is implicit with $\alpha=1$ only if
an exact algorithm exists.

\begin{theorem}
  Approximating the minimum \ran subject to a connectivity constraint is NP-Hard.
\end{theorem}

\pf Recall that a Hamiltonian cycle on $G$ is a tour (set of adjacent
edges or, equivalently, nodes in $G$) that visits each node exactly
once, except for the start node which is equal to the last node on the
tour. We claim the existance of a Hamiltonian cycle on a given graph
is equivalent to the feasibility of the minimum \ran with connectivity
on a related instance. Recall that an instance of Hamiltonian cycle consists of a graph, so
let such an instance be given with $G = (N,E)$. Now define the vector
$b \in \R^{|N|}$ by setting $b_i = 2$ for $i \in \{1,\ldots,|N|\}$ and
consider the resulting connected minimum \ran instance using the graph
$G$ and vector $b$. 

We first show that if the minimum \ran instance $(G,b)$ is feasible,
then there is a Hamiltonian cycle on $G$. Suppose there is a feasible
solution $H=(N,F)$, which means $F \subseteq E$, each node $u \in N$
has degree 2, and $H$ is connected. As each node has even degree and
$H$ is connected, there is an Eulerian cycle, $T$, on $H$. We claim
that $T$ is a Hamiltonian cycle on $G$, which means each node is
visited exactly once by $T$ except the start node which is visited
exactly twice. Choose $u \in N$ and note that two edges are adjacent to $u$. Then, because $T$ traverses every arc, the node $u$ is visited. Denote the start node of $T$ by $s \in N$ and consider traversing
$T$ beginning at $s$. If the traversal visits a node $u \in N
\setminus \{s\}$ more then once then an edge was traversed into $u$, a
second distinct edge was traversed out of $u$, and a third distinct
edge was traversed into $u$, a contradiction as there are exactly two
distinct edges adjacent to $u$ in $H$. The same argument applies if
$s$ is visited more then once before the traversal is complete. So
each node is visited exactly once by $T$ except the start node, which
is returned to when the traversal is complete, i.e., $T$ is a
Hamiltonian cycle. Thus, if the instance $(G,b)$ is feasible, then $G$
possess a Hamiltonian cycle.

We now show that if there is a Hamiltonian cycle on $G$, then the
minimum \ran problem on $(G,b)$ is feasible. Consider a Hamiltonian
cycle, $C$ on $G$ and the subgraph induced by $C$. Such a subgraph is
connected as each node is visited. Also, each node has degree 2 as
each node $u \in N$ has one arc used to enter $u$ and exactly one
distinct arc $u$ to exit. Thus, if $G$ possess a Hamiltonian cycle,
then $(G,b)$ must be feasible to the given connected minimum \ran
instance.

Now suppose there were an $\alpha$-approximation algorithm to the
connected minimum \ran problem for some $\alpha \geq 1$. If the
algorithm returns a solution to the instance $(G,b)$ then $G$ posesses
a Hamiltonian cycle. If it does not, then the instance $(G,b)$ was not
feasibile and $G$ does not posess a Hamiltonian cycle. Note that the
argument does not rely on what value $\alpha$ is.  $\ddag$

The following corollaries are immediate.
\begin{corollary}
  Approximating $b$-matching when the graph induced by the matching must
  be connected is NP-Hard.
\end{corollary}

Because the proof showed that finding a feasible solution to the minimum \ran index is NP-Hard, we immediately can state the following corollary.
\begin{corollary}
  Approximating the maximum \ran subject to a connectivity constraint is NP-Hard.
\end{corollary}

Note that we have not shown what the complexity is when the input
graph $G$ is the complete graph. We leave this as an open problem.

\subsection{Example Transformation} \label{subset:re} Given the degree
sequence $d=(3,2,2,2,2,1)$, what is a graph realization with the minimum \rann?
We let nodes $v_1,v_2,v_3,v_4,v_5, v_6 \in N$ with $b=(3,2,2,2,2,1)$. Now we can form the complete graph $G$,
with weights corresponding to $b_i\cdot b_j$ for every node $v_{i},v_{j} \in
N$.
\begin{center}
\begin{tikzpicture}[node distance=30mm]
        \tikzstyle{nnode}=%
        [%
          fill=blue!50!black!20,%
          draw=blue!50!black,%
          minimum size=8mm,%
          circle,%
          thick,%
          inner sep=0pt%
        ]
         
        \node[nnode] (a)              {$3$};
        \node (la) [xshift = -.7 cm]{\textcolor{red}{$v_1$}};
        \node[nnode]  (b) [right of=a] {$2$};
        \node (lb) [right of= a, xshift = .7 cm]{\textcolor{red}{$v_2$}};
        \node [nnode] (c) [below of=a, xshift= -1.5 cm] {$2$};
        \node (lc) [below of=a, xshift = -2.2 cm]{\textcolor{red}{$v_3$}};
         \node (lg) [below of=a, xshift = -3.2 cm]{$G:$};
        \node [nnode] (d) [below of=b, xshift=1.5 cm] {$2$};
        \node (ld) [below of=b, xshift = 2.2 cm]{\textcolor{red}{$v_4$}};
        \node (le) [below of=c, xshift = .7 cm]{\textcolor{red}{$v_5$}};
        \node[nnode]  (e) [below of=c, xshift = 1.5 cm] {$2$};
        \node[nnode]  (f) [below of=d, xshift = -1.5 cm] {$1$};
        \node (lf) [below of=d, xshift = -.7 cm]{\textcolor{red}{$v_6$}};

         \path [thick][every node/.style={font=\sffamily\small}]    
         		     %outside edges
        		      (a) edge node[above] {6} (b) 
		      (b) edge node[above, xshift=.25 cm] {4} (d) 
		      (d) edge node[below, xshift=.25 cm] {2} (f) 
		      (e) edge node[below] {2} (f) 
		      (e) edge node[below, xshift=-.25 cm] {4} (c) 
		      (c) edge node[above, xshift=-.25 cm] {6} (a)
		      % rest of a edges
		      (a) edge node[above, xshift=-1.25cm, yshift=.75cm] {6} (d) 
		      (a) edge node[xshift=-.8cm, yshift=2cm] {3} (f) 
		      (a) edge node[xshift=.2cm, yshift=1.8cm ]{6} (e)
		      %rest of b edges
		      (b) edge node[xshift=-.2cm, yshift=1.75cm] {2} (f) 
		      (b) edge node[xshift=.8cm, yshift=2.1cm] {4} (e)
		      (b) edge node[xshift=1.25cm, yshift=1.1cm] {4} (c) 
		      %rest of c edges
		      (c) edge node[xshift=-2cm, yshift=.25cm] {4} (d) 
		      (c) edge node[xshift=-1.25cm, yshift=1.1cm] {2} (f)
		      %rest of d and e edges
		      (d) edge node[xshift=1.25cm, yshift=1.1cm] {4} (e);

\end{tikzpicture} \\ 
\vskip .5 true in 
\end{center} 
  $$v_1~~ v_2~~ v_3 ~~v_4 ~~v_5 ~~v_6 $$
$$  \begin{bmatrix}
    0 & 6 & 6 & 6 & 6 & 3 \\ 
    6 & 0 & 4 & 4 & 4 & 2 \\ 
    6 & 4 & 0 & 4 & 4 & 2 \\ 
    6 & 4 & 4 & 0 & 4 & 2 \\ 
    6 & 4 & 4 & 4 & 0 & 2 \\ 
    3 & 2 & 2 & 2 & 2 & 0 \\ 
  \end{bmatrix}$$
Now we solve the minimum weight perfect $b$-matching for $G$ and obtain $G'$:
\begin{center}
\begin{tikzpicture}[node distance=30mm]
        \tikzstyle{nnode}=%
        [%
          fill=green!50!black!20,%
          draw=green!50!black,%
          minimum size=8mm,%
          circle,%
          thick,%
          inner sep=0pt%
        ]

        \node[nnode] (a)              {$3$};
        \node (la) [xshift = -.7 cm]{\textcolor{red}{$v_1$}};
        \node[nnode]  (b) [right of=a] {$2$};
        \node (lb) [right of= a, xshift = .7 cm]{\textcolor{red}{$v_2$}};
        \node [nnode] (c) [below of=a, xshift= -1.5 cm] {$2$};
        \node (lc) [below of=a, xshift = -2.2 cm]{\textcolor{red}{$v_3$}};
           \node (lg) [below of=a, xshift = -3.2 cm]{$G':$};
        \node [nnode] (d) [below of=b, xshift=1.5 cm] {$2$};
        \node (ld) [below of=b, xshift = 2.2 cm]{\textcolor{red}{$v_4$}};
        \node (le) [below of=c, xshift = .7 cm]{\textcolor{red}{$v_5$}};
        \node[nnode]  (e) [below of=c, xshift = 1.5 cm] {$2$};
        \node[nnode]  (f) [below of=d, xshift = -1.5 cm] {$1$};
        \node (lf) [below of=d, xshift = -.7 cm]{\textcolor{red}{$v_6$}};
        
         \path [thick][every node/.style={font=\sffamily\small}]    
         		     %outside edges
        		      (a) edge node[above] {6} (b) 
		      (c) edge node[above, xshift=-.25 cm] {6} (a)
		      % rest of a edges
		      (a) edge node[xshift=-.7cm, yshift=2cm] {3} (f) 
		      %rest of b edges
		      (b) edge node[xshift=.8cm, yshift=2.1cm] {4} (e)
		      %rest of c edges
		      (c) edge node[xshift=-2cm, yshift=.25cm] {4} (d) 
		      %rest of d and e edges
		      (d) edge node[xshift=1.25cm, yshift=1.1cm] {4} (e);

\end{tikzpicture} \\ 
\vskip .5 true in 
\end{center} 
$$  \begin{bmatrix}
    0 & 1 & 1 & 0 & 0 & 1 \\ 
    1 & 0 & 0 & 0 & 1 & 0 \\ 
    1 & 0 & 0 & 1& 0 & 0 \\ 
    0 & 0 & 1 & 0 & 1 & 0 \\ 
    0 & 1& 0 & 1 & 0 & 0 \\ 
    1 & 0 & 0 & 0 & 0 & 0 \\ 
  \end{bmatrix}$$
$$R(G') = 6+6+4+4+4+3 = 27 $$
$G'$ is a solution for the minimum weight perfect $b$-matching. 
The sum of the weights is the minimum \ran and the unweighted adjacency 
matrix is the corresponding graph realization. Note that there are other 
solutions to the matching that will produce the minimum \ran and a different 
realization. That is, the solution is not unique.

\section {Algorithms}

Our primary goal is to devise an algorithm to solve the minimum \ran
problem.  An algorithm that is useful when creating graphs with a
specified degree sequence is the Havel-Hakimi algorithm.  The
Havel-Hakimi algorithm can be used to check if a degree sequence is
graphic and to find a realization of that sequence.

\subsection{Havel-Hakimi Algorithm}
\label{sec:havel}
The Havel-Hakimi algorithm is useful when we have a non-negative integer sequence and we want to know if it is graphic, and if so, what is a realization of the sequence.  Its algorithm is given below:\\

\def\norm#1{\left\|#1\right\|}
\newcommand{\fprod}[1]{\ensuremath{\left\langle #1 \right\rangle}}
\begin{center}
   \fbox{
     \begin{minipage}{.95\linewidth}
       \begin{codebox}
         {\Procname{\underline{\proc{Havel-Hakimi}} \cite{havel1955remark,hakimi1962realizability}}}
         \zi {\bf Inputs: } $d$, a non-negative integer sequence\\
         \zi {\bf Outputs: } $A$, adjacency matrix of a graph realization of $d$ (if graphic) \\
         \zi Initialize a $|d| \times |d|$ adjacency matrix $A$ so that $A_{ij}=A_{ji} = 0$.
         \zi \While $d$ is not the 0 sequence
         \zi \Do
         \zi Pick a random index $i$ of $d$.
         \zi Subtract $1$ from the $d_i$ nodes with largest degrees (not including $d_i$) and set $d_i$ to 0.
         %\zi \ \ \ largest degree) until $x=0$.
         \zi For the $d_i$ nodes of largest degrees in the previous step (call them $v_j$), set 
         \zi \ \ \ $A_{ij} = A_{ji} = 1$.
	  \zi \If $|d| == 1$
	  \zi \ \ \ \ \ $d$ is not graphic.
	  \zi \ \ \ \ \ \Return null.
         \zi \End
         \zi \Return $A$
         \label{alg:HH}
       \end{codebox}
     \end{minipage}
   }
   \vspace*{+.1in}
\end{center}

\subsection{Two-Switches and the Metagraph} 
\label{sec:twoswitch}
One way to generate a collection of realizations for a degree sequence is to move 
from one realization to another by doing a two-switch, with an example as follows:
\begin{center} 
\begin{tikzpicture}[node distance=15mm]
        \tikzstyle{every node}=%
        [%
          fill=green!50!black!20,%
          draw=green!50!black,%
          minimum size=8mm,%
          circle,%
          thick,%
          inner sep=0pt%
        ]

        \node (A)              {$a$};
        \node (C) [right of=A] {$c$};
        \node (B) [below of=A] {$b$};
        \node (D) [below of=C] {$d$};

        \path [thick](A) edge (B)
                           (C) edge (D);

      \end{tikzpicture}
~~~~~~~~~~~~~
\begin{tikzpicture}[node distance=15mm]
        \tikzstyle{every node}=%
        [%
          fill=green!50!black!20,%
          draw=green!50!black,%
          minimum size=8mm,%
          circle,%
          thick,%
          inner sep=0pt%
        ]

        \node (A)              {$a$};
        \node (C) [right of=A] {$c$};
        \node (B) [below of=A] {$b$};
        \node (D) [below of=C] {$d$};

        \path [thick](A) edge (B)
                           (C) edge (D);
        \path [dashed](A) edge (D)
                         (C) edge (B);

      \end{tikzpicture}
~~~~~~~~~~~~~
\begin{tikzpicture}[node distance=15mm]
        \tikzstyle{every node}=%
        [%
          fill=green!50!black!20,%
          draw=green!50!black,%
          minimum size=8mm,%
          circle,%
          thick,%
          inner sep=0pt%
        ]

        \node (A)              {$a$};
        \node (C) [right of=A] {$c$};
        \node (B) [below of=A] {$b$};
        \node (D) [below of=C] {$d$};

        \path [thick](A) edge (D)
                         (C) edge (B);

      \end{tikzpicture}
\end{center}
When doing a two-switch, we examine two edges, $(a,b),(c,d)\in E$.  If $(a,d)\notin E$ and $(b,c)\notin E$ then we can remove edges $(a,b)$ and $(c,d)$ and create 
edges $(a,d)$ and $(b,c)$. This is not a unique move, since we could also use $(a,c)$ and $(b,d)$ if $(a,c)\notin E$ and $(b,d) 
\notin E$. Two-switching is an easy way to obtain a different graph with the same degree sequence after a graph is created using 
the Havel-Hakimi algorithm.  The two-switch pseudocode is given below:\\
\begin{center}
   \fbox{
     \begin{minipage}{.85\linewidth}
       \begin{codebox}
         {\Procname{\underline{\proc{Two-Switch}}}}
         \zi {\bf Inputs: }$G$, a graph with degree sequence $d$.
         \zi {\bf Outputs: }$G^{\prime}$, a new graph with degree sequence $d$.
         \zi Let $G^{\prime} := G$.
         \zi Pick a random edge $(a,b)\in E(G^{\prime})$.
         \zi Find a node $c\in N(G^{\prime})$ that is not connected to $b$.
         \zi Find a node $d\in N(G^{\prime})$ that is connected to $c$, but not to $a$.
         \zi Remove edges $(a,b),(c,d)$ from $E(G^\prime)$ and add $(a,d),(c,b)$ to $E(G^\prime)$.
         \zi \Return $G^{\prime}$
         \label{alg:tw}
       \end{codebox}
     \end{minipage}
   }
   \vspace*{+.1in}
\end{center}
We can construct a {\it metagraph} of a degree sequence, where
the metagraph is an undirected graph with each node representing a graph realization of a degree sequence and each edge representing a two-switch.  The following theorem shows that the metagraph is always a connected graph.
\begin{theorem} (Ryser's Theorem \cite{ryser1957combinatorial}.)
Given graphs $G$ and $G'$ such that $d(G)=d(G')$, there exists a sequence of two-switches going from $G$ to $G'$.
\end{theorem}

\subsection{Heuristic for Disconnected Realizations}
\label{sec:heuristic}

Since finding a connected graph realization with minimum \ran is NP-Hard,
we present a heuristic using two-switches to connect disconnected realizations. The heuristic sequentially performs a two-switch 
between pairs of connected components until all the components are connected:
\begin{center}
   \vspace*{+.1in}
   \fbox{
     \begin{minipage}{.85\linewidth}
       \begin{codebox}
         {\Procname{\underline{\proc{Two-switch Heuristic}}}}
         \zi {\bf Inputs: }$G$, a disconnected graph with degree sequence $d$
         \zi {\bf Outputs: }$G^{\prime}$, a connected graph with degree sequence $d$
         \zi
         \zi Let $G^{\prime} := G$.
         \zi \While the number of connected components in $G^{\prime}$ is $\geq 2$
         \zi \Do a two switch with two components to connect them using two randomly
         \zi chosen edges from each component
         \zi \End
         \zi \Return $G^{\prime}$
       \end{codebox}
     \end{minipage}
   }
   \vspace*{+.1in}
\end{center}
 
Note that the method to connect the disconnected realizations may not produce graphs 
with the best structure since there is only 1 edge connecting one component to another. 
Also note that we do not need to check whether the randomly chosen edges are adjacent or 
not since they are in separate connected components. 

\section{Solving the Minimum \ran Problem}

In this section, we focus on the case where the input graph $G$ is the
complete graph.
To solve the minimum \ran problem we used code that solves a minimum weight perfect $b$-matching problem. 
The code used is for generalized matching problems and was written by Vlad Schogolev, Bert Huang, and Stuart Andrews. 
Their code uses the GOBLIN graph library (http://goblin2.sourceforge.net/). Huang's paper on loopy belief propagation for bipartite maximum weight 
$b$-matching uses this code \cite{huang2007loopy}. 
The code solves a \textit{maximum} weight perfect $b$-matching problem given the weight matrix, $b$ vector and the exact solution algorithm choice. 
After transforming our minimum \ran problem instance into a minimum weight perfect 
$b$-matching instance, the $b$ vector will correspond to the degree constraints and the weight matrix to the possible degree products (see Section 
\ref{subset:re}). However, since the code solves the $maximum$ matching we transform our problem so that solving for the maximum yields the solution for the minimum. \\

Given a weight matrix $H$ we transform these weights into a matrix $H_2$ such that the maximum matching using $H_2$ will yield the same 
solution as the minimum matching using $H$. To do this we take a matrix $M$ with 1s in all positions except for the diagonal which has 0s. We 
then multiply every entry by one more than the maximum entry of $H$. $H$ is then subtracted from $M$ yielding $H_2$. 
The following pseudocode implements this algorithm and will solve the minimum \ran problem for a given degree sequence:
\begin{center}
   \vspace*{+.1in}
   \fbox{
     \begin{minipage}{.85\linewidth}
       \begin{codebox}
         {\Procname{\underline{\proc{Algorithm to solve minimum \ran with $b$-matching }}}}
         \zi {\bf Inputs: }$A$, an adjacency matrix with degree sequence $d$.
         \zi {\bf Outputs: }$A^{\prime}$, the new adjacency matrix with degree sequence $d$ 
         \zi  ~~~~~~~~~~~~~~~and minimized Randic index $r$.
         \zi
         \zi Create weight matrix $H$ of degree products.
         \zi Let $b := d$.
         \zi Let $H_2 := (1+\max H)\cdot M - H$, where $M$ is the adjacency matrix of all 1s except 
         \zi \ \ \ along the diagonal.
         \zi Call b-match solver with $b$ and $H_2$ to get adjacency matrix $A^{\prime}$ of optimal solution
         \zi Calculate $r=R(A^{\prime})$
         \zi \Return $A^{\prime}$ and $r$
       \end{codebox}
     \end{minipage}
   }
   \vspace*{+.1in}
\end{center}
This algorithm returns the minimum \ran of a graph and a realization. 
We know that the $b$-matching code runs in polynomial time (\cite{schrijver2003combinatorial}).
It is easy to see that the transformation steps are done in 
polynomial time as well. We used three types of randomly generated graphs to test the algorithm performance: 
\ern, geometric and scale-free.  We limited our computational experiments to degree sequences for which connected realizations
were known to exist. The \ran before and after the optimization was recorded.   After the optimization we check 
if the graph realization with the minimum \ran is connected. 
We generated graphs with 25, 50 and 100 nodes. In addition, 100 of each graph type and size were generated. 
Tables \ref{tab:25}, \ref{tab:50} and \ref{tab:100}  present  results from the runs. Note that the number of graphs connected 
after the run plus the number of graphs disconnected plus the number of graphs with no connected realizations is 100 for each graph type. 
\begin{table}[h!]
  \centering 
  \begin{tabular}{c|c|c|c}
%    \hline
{\bf    Graph type} & {\bf connected} & {\bf disconnected} & {\bf no connected realizations}\\ 
    \hline
    \er & 67 & 1 & 32\\ 
%\hline
    Geometric & 61 & 2 & 37\\ 
%\hline
    Scale-Free & 93 & 7 & 0\\ 
    \hline
  \end{tabular}
  \caption{25 node graphs}
  \label{tab:25}
\end{table}
\begin{table}[h!]
  \centering
  \begin{tabular}{c|c|c|c}
%    \hline
{\bf    Graph type}&{\bf connected} & {\bf disconnected} & {\bf no connected realizations}\\ 
    \hline
    \er & 50 & 5 & 45\\ 
%\hline
    Geometric & 57 & 3 & 40\\ 
%\hline
    Scale-Free & 85 & 15 & 0\\ 
    \hline
  \end{tabular}
  \caption{50 node graphs}
  \label{tab:50}
\end{table}
\begin{table}[h!]
  \centering
  \begin{tabular}{c|c|c|c}
%    \hline
{\bf    Graph type}&{\bf connected} &{\bf disconnected} &{\bf no connected realizations}\\ 
    \hline
    \er & 16 & 2 & 82\\ 
%\hline
    Geometric & 30 & 6 & 64\\ 
%\hline
    Scale-Free & 91 & 8 & 1\\ 
    \hline
  \end{tabular}
  \caption{100 node graphs}
  \label{tab:100}
\end{table}
The MATLAB functions used to generate the geometric and scale-free graphs are from CONTEST: 
A Controllable Test Matrix Toolbox for MATLAB \cite{taylor2009contest}. 
In addition, the necessary and sufficient conditions for a non-negative integer sequence $\{a_i\}$ to be 
realizable as the degrees of the nodes of a connected graph are that $a_i \neq 0$ for all $i$ 
and the sum of the integers $a_i$ is even and not less than $2(n-1).$ This condition was used to 
discard graphs with a degree sequence that had no connected realizations \cite{chen1997}.\\

In general from our runs, the realization of the minimum \ran was connected. There are minimum \ran 
graph realizations that are disconnected and we do not know if there are other realizations with this 
\ran that are connected since the $b$-matching solver only produces one solution. But there were often 
a large proportion of graphs that had no connected realization at all. This largely depends on parameters 
chosen for the randomly generated graphs.   If the random graph produced has most nodes with large degrees 
then it is unlikely that any graph realization would be disconnected. We were interested in generating graphs 
that have both connected and disconnected realizations and investigating whether the realization generated with 
the minimum \ran was connected or not.\\

The following parameters were chosen after extensive experimentation so that the number of instances with no
connected realizations was small.
For the \er graphs we used an average degree per node of 4.25. The corresponding $p$ values used were calculated 
using $p = \frac{4.25}{n}$ where $n$ is the number of nodes in the graph. Thus $p=.17$ for $n=25$, $p=.085$ for 
$n=50$, and $p=.043$ for $n=100$. For the geometric graphs we used an average degree per node of 6. The radii 
were calculated using $r=\sqrt{\frac{6}{\pi n}}$. Our corresponding radii were $r=.276$ for $n=25$, $r=.195$ 
for $n=50$, and $r=.138$ for $n=100$. We used scale-free graphs with a minimum node degree of 2.\\

\begin{figure}[!ht]
   \centering
   \resizebox{.6\textwidth}{!}{\includegraphics{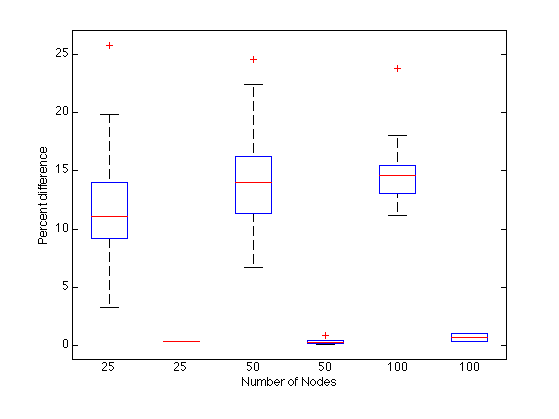}}
      \caption{Comparing percent differences for \er graphs.}
   \label{fig:er3}
\end{figure}

\begin{figure}[!h]
   \centering
   \resizebox{.6\textwidth}{!}{\includegraphics{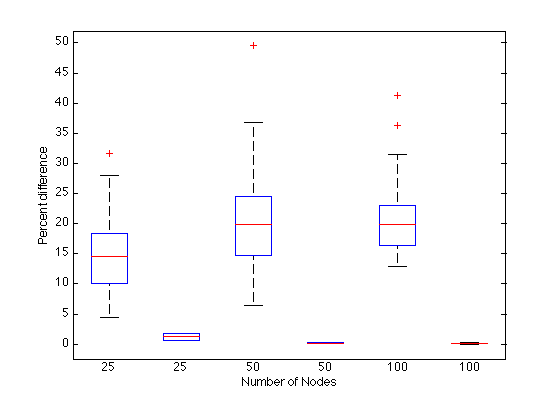}}
      \caption{Comparing percent differences for geometric graphs.}
   \label{fig:geo3}
\end{figure}

\begin{figure}[!htb]
   \centering
   \resizebox{.6\textwidth}{!}{\includegraphics{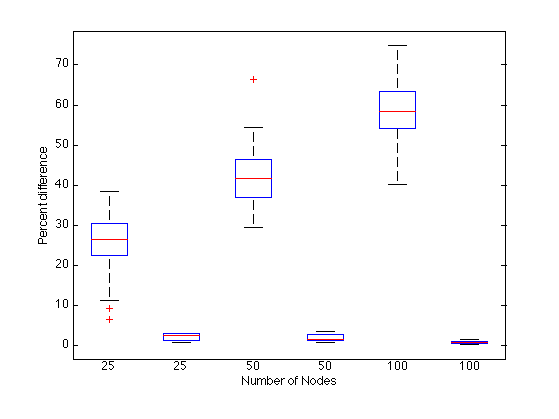}}
      \caption{Comparing percent differences for scale-free graphs.}
   \label{fig:pref3}
\end{figure}

The left box plots for each of 25, 50 and 100 node graphs in Figures \ref{fig:er3}, \ref{fig:geo3} and \ref{fig:pref3} show the percent difference between 
the graph's original \ran and the minimum \rann. The percent difference is calculated from $\frac{original - minimum} {minimum} \times 100.$ 
The right box plots for each of 25, 50 and 100 node graphs in Figures \ref{fig:er3}, \ref{fig:geo3}, and \ref{fig:pref3} show the percent difference between the minimum \ran and the \ran after the heuristic algorithm in Section~\ref{sec:heuristic} was applied. 
This percent 
difference is calculated with $\frac{after~heuristic - minimum} {minimum} \times 100.$ The number of graphs that used the heuristic 
depended on the number of optimal graph realizations that were disconnected. Note that this is a different number for each graph type 
and size. See Tables \ref{tab:25}, \ref{tab:50} and \ref{tab:100} for those numbers. 

\section{Applications}
\label{sec:apps}

We will now show results for two applications where the effects of network connectivity measures have shown or are hypothesized to play an important role:  neuronal synchronization and dementia.

\subsection{Neuronal synchronization}
\label{sec:netsync}

The first example is exploring the effect of assortativity on the synchronous firing of neurons in the preB\"{o}tzinger complex.  This collection of neurons is responsible for the control of respiration in mammals \cite{Feldman:2006p907}.  As neuronal networks are inherently directed, we must first define a directed \ran \cite{ZamoraLopez:2008p6831} and extend the $b$-matching algorithm to the directed case.  The extension is similar to \cite{liuNature2011} in which the maximum matching problem is extended to directed graphs for a network controllability problem.

\subsubsection{Extension to directed networks}

\begin{figure}
\centering
\includegraphics{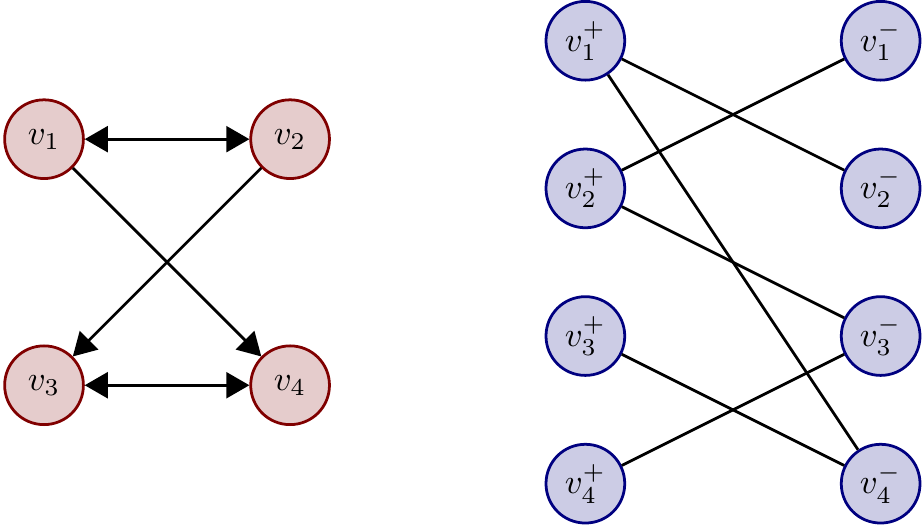}
\caption{{\bf Left}: Directed graph $\vec{G}=(N,E)$ with node set $N=\{v_{1},v_{2},v_{3},v_{4}\}$ and degree sequence $\{(2,1),(2,1),(1,2),(1,2)\}$.  {\bf Right}:  Bipartite representation $\vec{G}^{*} = (N^{*},E^{*})$, where $N^{*} = N^{\p}\cup N^{\m}$, $N^{\p} =  \{v_{1}^{\p},\ldots,v_{n}^{\p}\}$, $N^{\m} =  \{v_{1}^{\m},\ldots,v_{n}^{\m}\}$, and $e^{*} = (v_{i}^{\p},v_{j}^{\m}) \in E^{*}$ if and only if $e = (v_{i},v_{j}) \in E$.}
\label{fig:digraph_example}
\end{figure}

Consider the directed graph $\vec{G} = (N,E)$ with vertex set $N = \{v_{1},\ldots,v_{n}\}$ and edge (or arc) set $E = \{(v_{i},v_{j}) \, | \, v_{i} \rightarrow v_{j}\}$.  The degree sequence for $\vec{G}$ is a non-negative integer-pair sequence $d = \{(d^{\p}_{i},d^{\m}_{i}) \, | \, i=1,\ldots,n\}$, where we denote the out-degree and in-degree sequences by $d^{\p}$ and $d^{\m}$, respectively.  There are four different \ran-type measures \cite{ZamoraLopez:2008p6831} given by 
\begin{equation}
\label{eq:dirandic}
R^{pq}(G) = \sum_{(v_i,v_j)\in E}d^{p}_{i}d^{q}_{j},
\end{equation}
where $p,q\in\{-,+\}$.  This can be seen as a natural extension of $R(G)$ to the directed case.  In this context, we will now define an extension of perfect $b$-matching as follows.
\begin{definition}[Perfect $b$-Matching for Directed Graphs]
\label{def:dir_pbmatch}
For a directed graph $\vec{G} = (N,E)$ and positive integer-pair sequence $b = (b^{\p},b^{\m}) = \{(b^{\p}_i,b^{\m}_{i}) \in \Z^n \times Z^n \, | \, i=1\,\ldots,n\}$, a {\it perfect $b$-matching} is a subset of edges $M
\subseteq E$ such that for node $v_i \in N$, the out and in-degree of $v_i$ in the subgraph $(N,M)$ is $b^{\p}_i$ and $b^{\m}_{i}$, respectively.  
\end{definition}
\noindent In an analogous way as the undirected case, we can define a minimum weight perfect b-matching problem as follows.
\vskip 0.1in
{\noindent}\textbf{Minimum Weight Perfect b-Matching Problem.}  Given a positive integer-pair sequence $b$, a directed graph $\vec{G}=(N,E)$ and a set of edge weights $w:E \rightarrow \R$, find a perfect $b$-matching with minimum weight.
\vskip 0.1in
In order to use the existing algorithm for undirected graphs, we consider the equivalent bipartite form of $\vec{G}$ given by $\vec{G}^{*} = (N^{*},E^{*})$, where $N^{*} = N^{\p}\cup N^{\m}$, $N^{\p} =  \{v_{1}^{\p},\ldots,v_{n}^{\p}\}$, $N^{\m} =  \{v_{1}^{\m},\ldots,v_{n}^{\m}\}$, and $e^{*} = (v_{i}^{\p},v_{j}^{\m}) \in E^{*}$ if and only if $e = (v_{i},v_{j}) \in E$ (see Figure~\ref{fig:digraph_example} for an example).
We can define a perfect $b$-matching on the corresponding bipartite graph $\vec{G}^{*} = (N^{*},E^{*})$ as a subset of edges $M^{*}\subseteq E^{*}$ such that for node $v_i \in N^{\p}$ or $N^{\m}$, the degree of node $v_i$ in $(N^{*},M^{*})$ is $b^{\p}_{i}$ or $b^{\m}_{i}$, respectively.  This modified definition is a special case of the undirected version in Definition \ref{def:pbmatch} with node set $N = (v_{1}^{\p},\ldots,v_{n}^{\p},v_{1}^{\m},\ldots,v_{n}^{\m})$ and the positive vector $b = (b_{1}^{\p},\ldots,b_{n}^{\p},b_{1}^{\m},\ldots,b_{n}^{\m})$.  Thus, to find a minimum weight perfect $b$-matching for a directed graph $\vec{G}$, we do the following (see Figure~\ref{fig:bigraph}):
\begin{figure}
\centering
\includegraphics{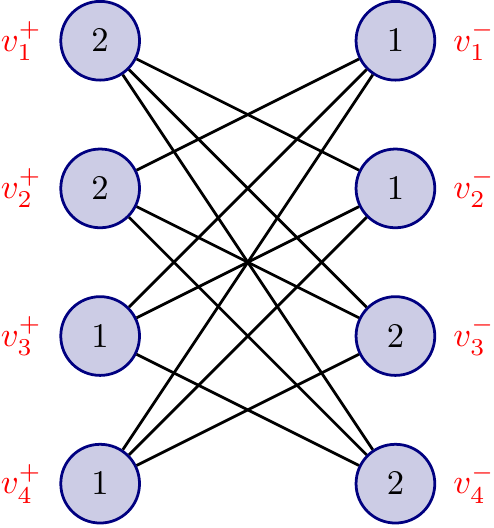}
\caption{Bipartite setup of the perfect $b$-matching problem for the directed graph given in Figure~\ref{fig:digraph_example}, with node set $N = (v_{1}^{\p},\ldots,v_{4}^{\p},v_{1}^{\m},\ldots,v_{4}^{\m})$ and the positive vector $b = (b_{1}^{\p},\ldots,b_{4}^{\p},b_{1}^{\m},\ldots,b_{4}^{\m}) = (2,2,1,1,1,1,2,2)$.  For minimum weight perfect $b$-matching, we would set $w_{ij} = d_{i}^{p}d_{j}^{q}$, where $p,q\in\{\p,\m\}$.}
\label{fig:bigraph}
\end{figure}
\begin{itemize}
\item Let $B$ be the complete bipartite graph $K_{n,n}$ minus the edges $\{(v_{i}^{\p},v_{i}^{\m}) \, | \, i=1,\ldots,n\}$.
\item Let $b^{\p} = d^{\p}$ and $b^{\m} = d^{\m}$.
\item For $p,q\in\{\p,\m\}$, let edge weights $w_{ij} = d_{i}^{p}d_{j}^{q}$.
\end{itemize}

\subsubsection{Neuronal networks}

\begin{figure}[!h]
\centering
\includegraphics[width=5in]{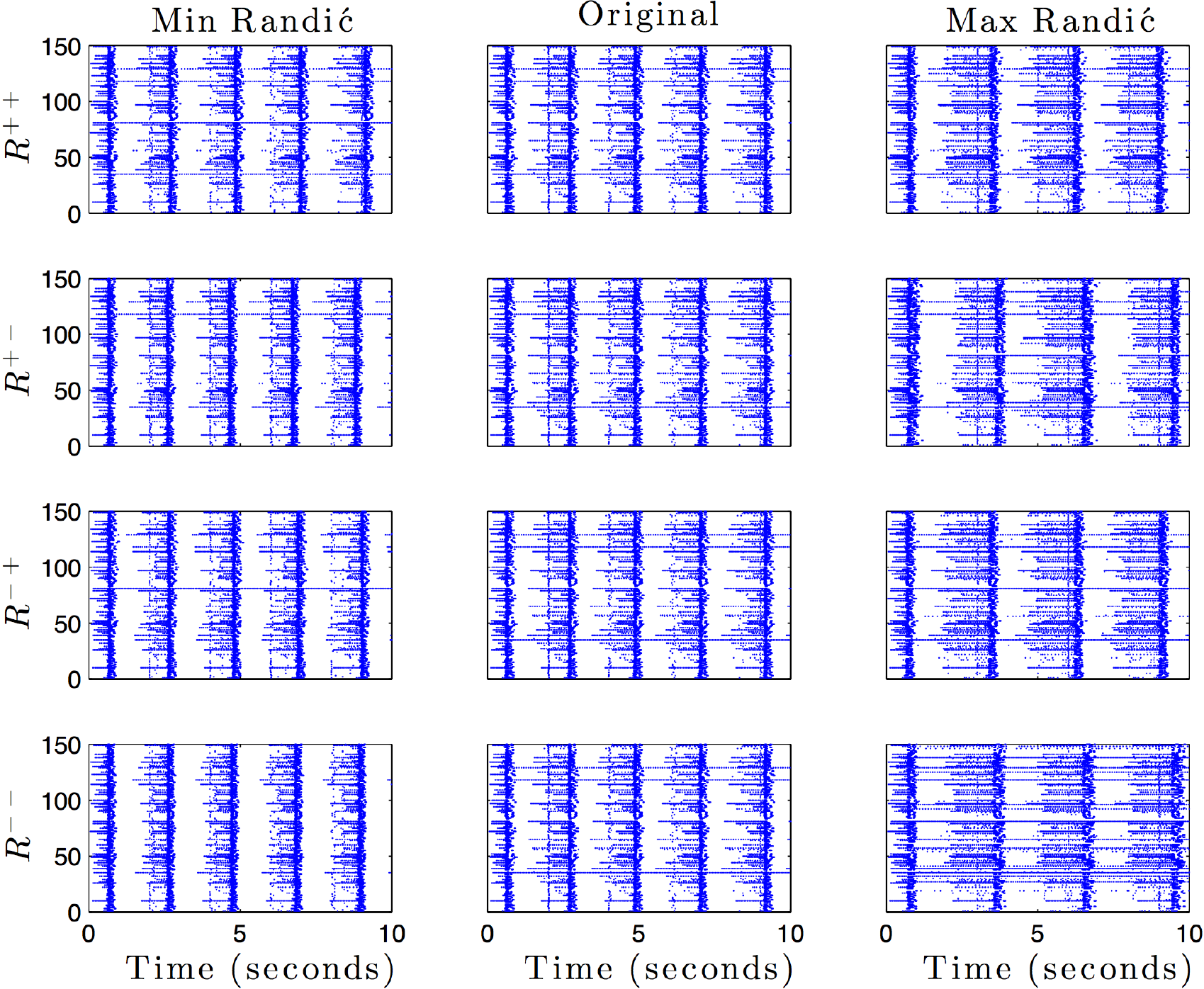}
\caption{Raster plots for simulations of 150 preB\"{o}tzinger complex neurons connected in a 3D geometric network (middle column, all rows the same), as well as network realizations with the minimum (left column) and maximum (right column) \ran $R^{pq}$, where $pq = +\,+$, $+\,-$, $-\,+$, and $-\,-$ in rows 1, 2, 3 and 4, respectively.}
\label{fig:neuron5}
\end{figure}

\begin{figure}[!h]
\centering
\includegraphics[width=5in]{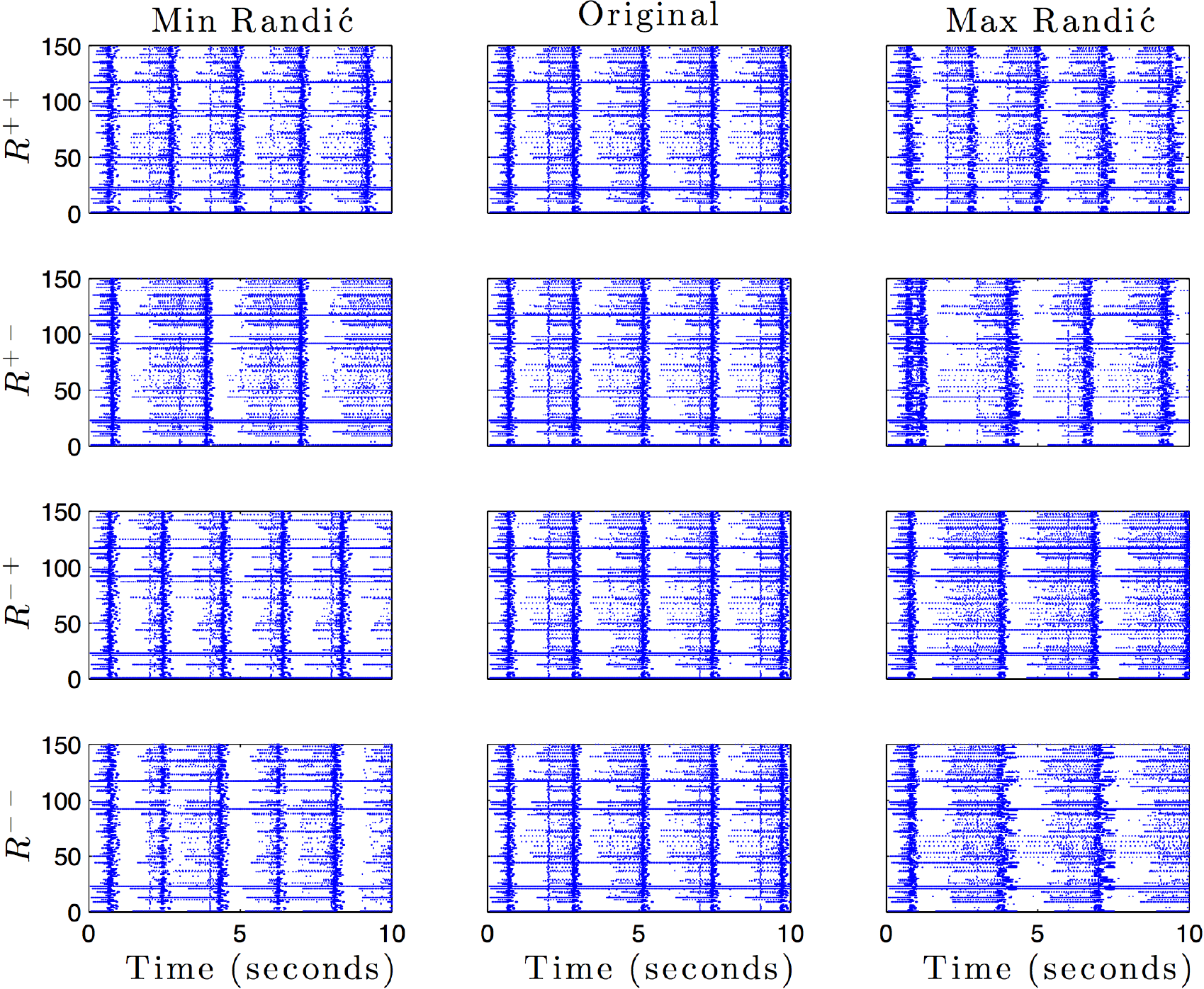}
\caption{Setup same as in Figure~\ref{fig:neuron5} for a different 3D geometric network.}
\label{fig:neuron20}
\end{figure}

Using NeuronetExperimenter \cite{NeuronetExperimenter:ra0SJb3l}, we simulated 150 rhythmogenic neurons in the preB\"{o}tzinger complex using the Rubin-Hayes neuron model \cite{Rubin:2009ie} (\url{http://senselab.med.yale.edu/modeldb/ShowModel.asp?model=125649}).  It is unknown what degree distribution and network connectivity arise in neuronal networks.  However, it seems reasonable to expect that neurons closer to each other are more likely to be connected.  Hence, we have chosen to model these networks with 3D geometric (directed) graphs.  Raster plots of simulation results for two 3D geometric networks are displayed in Figures~\ref{fig:neuron5} and \ref{fig:neuron20} (center columns), along with simulation results from realizations with the minimum and maximum directed Randi\'{c} indices in \eqref{eq:dirandic}.  These results display an overall tendency for the realization with minimum and maximum \ran, regardless of type, to lead to a faster and slower breathing rhythm, respectively.  This is not always the case, however, as can be seen with $R^{+\,-}$ in Figure~\ref{fig:neuron20}.  There is also a tendency for more and less synchronous firing from realizations with respective minimum and maximum \ran.  These results warrant a more thorough quantitative analysis, which is beyond the scope of this work.

\vfill\eject

\subsection{Normalizing the \ran index in cortical thickness networks}
\label{sec:ct}
Recently, researchers have begun using graph measures of connectivity to investigate the difference between structural magnetic resonance images taken from healthy individuals and individuals diagnosed with dementia (for example, \cite{Bai12,Yao10,HeChen09,HeChen08}). Networks in these studies are formed by calculating correlations between the cortical thicknesses of different brain regions, a technique based on the correlation between cortical thicknesss loss and dementia \cite{Fjell06,Gross12}. This method of network creation has been speculated to provide more insight on the functional relationships between brain regions \cite{Gong12}. Details of how the networks are formed can be found in \cite{HeChen08,SoldanMRJP14} and we present one such method below. After networks for each subject population (e.g., normal, subjects with dementia) are formed, network measures (e.g., \ran) are calculated and evaluated for significance.  

A particular challenge to using comparisons between the network measures is complicated by the varying number of edges in each of the networks. In particular, some kind of normalization is required for some of the measures used.  In order to normalize, we propose the following scheme: 
\begin{enumerate}
	\item For each network, $N$, determine the underlying degree distribution, $b$.
    \item Calculate $U_b$, the maximum \ran index for $b$, respectively.
    \item Use $U_b$ to normalize the \ran for $N$.
\end{enumerate}
We use data from the Alzheimer's Disease Neuroimaging Initiative (ADNI) and further analyzed by FreeSurfer, a technique developed by Fischl and Dale \cite{FischlD00} to measure cortical thicknesses. When applied to the ADNI data, Fischl and Dale discretized the cortical layer into 68 different regions. Soldan, et al.~\cite{SoldanMRJP14} used the following steps to generate networks:
\begin{enumerate}
\item Subset the population into categories based on whether they were diagnosed as normal (NORMAL), diagnosed with mild cognitive impairment for three or more years without disease progression (MCI), diagnosed with mild cognitive impairment and then progressed to Alzheimer's Disease within three years (MCIAD), or diagnosed with Alzheimer's Disease.
\item Within each population, use regression to control for subject age, gender, education level, and interaction effects between age and gender.
\item Use either partial or Pearson's correlations to calculate coefficients and $p$-values between each of the 68 regions.

\label{li:corr}
\item Use False Discover Rate \cite{Benjamini95} calculations to determine significant correlation coefficients with an error rate of 5\%.
\item Use one of the following schemes to determine edge weights: \label{li:weights}
	\begin{enumerate}
    \item Use no weights: simply include edges or not.
    \item Use the absolute value of the correlation coefficients as the edge weight. \label{li:abs}
    \item Use the product of the normalized cortical thicknesses (i.e., the cortical thicknesses divided by the maximimum)
    between the two regions connected by the significant edge. \label{li:prod}
    \item Use both \ref{li:abs} and \ref{li:prod}.
    \end{enumerate}
\end{enumerate}
Given the four diagnostic categories, the two possibilities in step \ref{li:corr}, and the four possibilities in step \ref{li:weights}, there are a total of 8 different networks. After calculating the associated degree distributions, we calculated the maximium \ran. As an example, we display the results for one network in Table \ref{tab:optimizedMRI}. The results for all eight networks can be seen in Appendix \ref{app:maxct}

\begin{table}
  \centering
  \begin{tabular}{l|c|c|c}
    \hline	Group	&	Original	&	Maximum ($U_b$)	& 	Normalized \\	\hline
AD	&	172910.5	&	182773.6	&	$94.6\%$	\\
MCI-AD	&	182098.8	&	192005.4	&	$94.84\%$	\\
MCI	&	146719	&	155715.5	&	$94.22\%$	\\
Normal	&	231957.8	&	244387.6	&	$94.91\%$	\\	\hline
  \end{tabular}
  \caption{Original, maximum, and normalized \ran. CT and absolute correlation weighting used with Pierson's correlations.}
  \label{tab:optimizedMRI}
\end{table}

By normalizing the \ran, comparisons between the groups can more accurately determine whether the assortativity was due to the actual network topology versus other features, such as the total number of edges. Currently, studies instead arbitrarily delete edges in networks, which effectively ignore different brain regions in order to compare non-normalized graph measures \cite{Yao10,HeChen09,HeChen08}. Normalizing by dividing out by the optimized metric allows for comparisons without ignoring network features. To actually use the optimization in such a statistical study, we note that significance testing (e.g., permutation testing to determine whether the comparisons are valid) would be required.

\section{Conclusions and Future Work}
We have shown that the minimum \ran problem can be solved in polynomial time.
With use of available $b$-matching code we have developed an algorithm that produces 
a graph realization with the minimum \ran for a given degree sequence.  Not all optimal 
solutions are connected.  A two-switch heuristic was developed to connect disconnected 
optimal solutions.  Although the graph structure of these new connected graphs is fragile, 
the \ran changed relatively little. \\

From our experiments, when generating \ern, geometric and scale-free graphs, the realization with 
the minimum \ran is generally connected.  Undoubtably, this result is influenced by the parameters chosen  
for the randomly generated graphs, but many graphs that will have a disconnected minimum \ran realization 
have no connected realizations at all. \\

There are a number of future topics to explore.  We want to develop a
better way to connect graphs using the two-switch heuristic so that
the structure of the graph is less fragile.  Further experiments with
the parameters of randomly generated graphs are needed to understand
the conditions under which the number of graphs that are disconnected
or have no connected realizations changes.  We also are interested in determining the complexity
of the connected \ran problem when the input graph is not the complete
graph.

\bibliographystyle{plain}
\bibliography{main}

\appendix

\section{Maximized assortativity in cortical thickness networks}
\label{app:maxct}

\newcommand{\rwt}{{\bf Abs. corr. wtd.}\xspace}
\newcommand{\ctwt}{{\bf CT wtd.}\xspace}
\newcommand{\partials}{{\bf Partials}\xspace}

We present the computational results from the different cortical thickness (CT) networks. Recall that there were eight different networks that were created based on the following three techniques:
\begin{enumerate}
\item \rwt: Edges weighted by the absolute value of the correlation;
\item \ctwt: Edges weighted by the product of the CT at each node; and
\item \partials: {\it Partial correlations} used to control for the effects of nodes not adjacent to each edge.
\end{enumerate}
In the tables below, we use \rwt, \ctwt, and \partials to denote that the associated technique was involved in the data creation. Thus, the default method for network creation would involve no weighting and ordinary Pierson correlations calculated between different brain regions. 

\begin{table}[b]
  \centering
  \begin{tabular}{l|c|c|c}
    \hline	Group	&	Original \ran	&	Maximum \ran	& 	Percentage	\\	\hline
AD	&	6327146	&	6327662	&	$99.99\%$	\\
MCI-AD	&	6004728	&	6005871	&	$99.98\%$	\\
MCI	&	4354002	&	4366974	&	$99.7\%$	\\
Normal	&	5730156	&	5732158	&	$99.97\%$	\\	\hline
  \end{tabular}
  \caption{Original and maximized \ran.}
  \label{tab:nnn}
\end{table}

\begin{table}
  \centering
  \begin{tabular}{l|c|c|c}
\hline	Group	&	Original \ran	&	Maximum \ran	& 	Percentage	\\	\hline
AD	&	210683	&	226302	&	$93.1\%$	\\
MCI-AD	&	82102	&	89411	&	$91.83\%$	\\
MCI	&	54245	&	59447	&	$91.25\%$	\\
Normal	&	11329	&	13183	&	$85.94\%$	\\	\hline
  \end{tabular}
  \caption{Original and maximized \ran. \partials used.}
  \label{tab:nny}
\end{table}

\begin{table}
  \centering
  \begin{tabular}{l|c|c|c}
\hline	Group	&	Original \ran	&	Maximum \ran	& 	Percentage	\\	\hline
AD	&	172910.5	&	182773.6	&	$94.6\%$	\\
MCI-AD	&	182098.8	&	192005.4	&	$94.84\%$	\\
MCI	&	146719	&	155715.5	&	$94.22\%$	\\
Normal	&	231957.8	&	244387.6	&	$94.91\%$	\\	\hline
  \end{tabular}
  \caption{Original and maximized \ran. \ctwt used.}
  \label{tab:nyn}
\end{table}

\begin{table}
  \centering
  \begin{tabular}{l|c|c|c}
\hline	Group	&	Original \ran	&	Maximum \ran	& 	Percentage	\\	\hline
AD	&	826.0913	&	1056.366	&	$78.2\%$	\\
MCI-AD	&	966.2843	&	1295.196	&	$74.61\%$	\\
MCI	&	615.2327	&	809.9756	&	$75.96\%$	\\
Normal	&	157.3287	&	231.8221	&	$67.87\%$	\\ \hline
  \end{tabular}
  \caption{Original and maximized \ran. \ctwt and \partials used.}
  \label{tab:nyy}
\end{table}

\begin{table}
  \centering
  \begin{tabular}{l|c|c|c}
\hline	Group	&	Original \ran	&	Maximum \ran	& 	Percentage	\\	\hline
AD	&	628898.9	&	653754.5	&	$96.2\%$	\\
MCI-AD	&	637520.3	&	664594.7	&	$95.93\%$	\\
MCI	&	305088.9	&	325998.4	&	$93.59\%$	\\
Normal	&	496892.5	&	526608.8	&	$94.36\%$	\\ \hline
  \end{tabular}
  \caption{Original and maximized \ran. \rwt used.}
  \label{tab:ynn}
\end{table}

\begin{table}
  \centering
  \begin{tabular}{l|c|c|c}
\hline	Group	&	Original \ran	&	Maximum \ran	& 	Percentage	\\	\hline
AD	&	1293.972	&	1651.687	&	$78.34\%$	\\
MCI-AD	&	1491.38	&	1908.273	&	$78.15\%$	\\
MCI	&	615.0955	&	809.4072	&	$75.99\%$	\\
Normal	&	129.3011	&	193.4662	&	$66.83\%$	\\	\hline
  \end{tabular}
  \caption{Original and maximized \ran. \rwt and \partials used.}
  \label{tab:yny}
\end{table}

\begin{table}
  \centering
  \begin{tabular}{l|c|c|c}
    \hline	Group	&	Original \ran	&	Maximum \ran	& 	Percentage	\\	\hline
AD	&	16741.64	&	18608.67	&	$89.97\%$	\\
MCI-AD	&	19398.13	&	21562.17	&	$89.96\%$	\\
MCI	&	10007.96	&	11262.25	&	$88.86\%$	\\
Normal	&	18109.99	&	20315.58	&	$89.14\%$	\\	\hline
  \end{tabular}
  \caption{Original and maximized \ran. \rwt and \ctwt used.}
  \label{tab:yyn}
\end{table}

\begin{table}
  \centering
  \begin{tabular}{l|c|c|c}
    \hline	Group	&	Original \ran	&	Maximum \ran	& 	Percentage	\\	\hline
AD	&	33.39597	&	54.8791	&	$60.85\%$	\\
MCI-AD	&	40.12877	&	65.14792	&	$61.6\%$	\\
MCI	&	26.36858	&	45.19809	&	$58.34\%$	\\
Normal	&	5.027669	&	10.57358	&	$47.55\%$	\\	\hline
  \end{tabular}
  \caption{Original and maximized \ran. \rwt, \ctwt, and \partials used.}
  \label{tab:yyy}
\end{table}

%%%%%%%%%%%%%%%%
\end{document}